\ifx\selectfont\undefined
\documentstyle[12pt]{article}
\else
\documentstyle[12pt,oldlfont]{article}
\fi

\newcount\hh
\newcount\mm
\mm=\time
\hh=\time
\divide\hh by 60
\divide\mm by 60
\multiply\mm by 60
\mm=-\mm
\advance\mm by \time
\def\hhmm{\number\hh:\ifnum\mm<10{}0\fi\number\mm}
\makeatletter
\@addtoreset{equation}{section}
\def\subsection{\@startsection{subsection}{2}{\z@}{-3.25ex plus-1ex
    minus-.2ex}{-.1em}{\reset@font\large\bf}}
\def\subsubsection{\@startsection{subsubsection}{3}{\z@}{-3.25ex plus
 -1ex minus-.2ex}{-.1em}{\reset@font\normalsize\bf}}
\makeatother
\title{Asymptotic $l_p$ spaces and bounded distortions}
\author{Vitali D. Milman \and  Nicole Tomczak-Jaegermann}
\date{}

\newtheorem{fact}{Fact}
\newtheorem{thm}{Theorem}
\newtheorem{prop}{Proposition}
\newtheorem{lemma}{Lemma}

\newbox\nrmbox
\setbox\nrmbox=\hbox{$\Vert$}
\def\nrmrule{\vrule height\ht\nrmbox depth1.2\dp\nrmbox}
\setbox\nrmbox=%
  \hbox{\kern0.15em\nrmrule\kern0.15em\nrmrule\kern0.15em\nrmrule\kern0.15em}
\newcommand{\Snorm}[1]%
  {\copy\nrmbox#1\copy\nrmbox\kern-0.03em\lower.4ex\hbox{}}
\newcommand{\rem}{\noindent{\bf Remark{\ \ }}}
\newcommand{\proof}{{\noindent\bf Proof{\ \ }}}
\newcommand{\qed}{\bigskip\hfill\(\Box\)}
\newcommand{\Rn}[1]{\mbox{{\it I\kern -0.25emR}$^{\,{#1}}$}}
\newcommand{\NN}{\mbox{{\it I\kern -0.25emN}}}
\newcommand{\QQ}{\mbox{\it Q\hspace*{-1.1ex}\rule{0.15ex}%
       {1.5ex}\hspace*{1.1ex}}}
\newcommand{\QQNP}{\QQ_{+}^{(\NN)}}
\newcommand{\QQN}{\QQ^{(\NN)}}
\newcommand{\RRN}{\Rn{(\NN)}}

\newcommand{\bt}{\beta}
\newcommand{\de}{\delta}
\newcommand{\ep}{\varepsilon}

\newcommand{\vt}{\theta}
\newcommand{\la}{\lambda}

\newcommand{\spn}{{\rm span\,}}
\newcommand{\dist}{{\rm dist\,}}
\newcommand{\dis}{d}
\newcommand{\supp}{{\rm supp\,}}

\newcommand{\codim}{{\rm codim\,}}
\newcommand{\conv}{{\rm conv\,}}

\newcommand{\cB}{{\cal B\,}}

\newcommand{\ie}{{\em i.e.,\/}}
\newcommand{\eg}{{\em e.g.,\/}}
\newcommand{\cf}{{\em cf.\/}}
\newcommand{\asm}{asymptotic}
\newcommand{\as}{{\em \asm\/}}
\newcommand{\aslym}{\asm}

\newcommand{\aim}{almost  isometric} 
\newcommand{\ai}{{\em \aim\/}}

\newcommand\As[3]{#1_#2^{\,(#3)}}
\begin{document}
\maketitle
\begin{abstract}
The new class of Banach spaces, so-called asymptotic
$l_p$ spaces, is  introduced and it is shown that
every Banach space with bounded distortions contains a subspace
from this class.
The proof is based on an investigation  of
certain functions, called enveloping functions,
which are intimately connected
with stabilization properties  of the norm.
\end{abstract}


\section{Introduction}

During the last year several problems of infinite-dimensional
Banach space theory, which remained  open for decades,
have been  finally solved.
Some new constructions
of Banach spaces
have been made  which, on one hand,
showed  limitations of the theory, but  on the other hand,
also showed how exciting an infinite-dimensional geometry
can be.
Let us mention few of them:
\begin{description}
\item[(i)]
  a space without unconditional basic sequence (Gowers--Maurey),
\item[(ii)]
  a space not isomorphic to any of  its hyperplanes (Gowers),
\item[(iii)]
  a space such that every bounded operator being a Fredholm operator
  (Gowers--Maurey).
\end{description}

The problems which were answered by these examples are
of a linear-to\-po\-lo\-gi\-cal nature.
Although  a thorough study of
this kind of properties
flourished  back in the 60s, methods developed
that time and later were not sufficient
to succesfully atack these problems. The solutions
given
last year
are  by-products
of a study in a different direction:
the  infinite-dimensional geometry of convex bodies,
that is,  the geometry of the unit
sphere of a Banach space.

In this introduction we would like to explain  the
geometry which led to the breakthrough described above;
and in the main body of the paper we would like to add
some information in this geometric direction.

\subsection{}
\label{spectrum}

In fact, the topic of studies which led to the recent development
takes its roots, in a large part, in  the local theory
of Banach spaces, in other words, in the asymptotic theory
of finite-dimensional normed spaces. Consider
the following question:

{\em Let $f(\cdot)$ be a uniformly continuous real valued function
on the unit sphere $S= S(X)= \{x\in X \mid  \|x\| =1\}$ of an
infinite-dimensional Banach space $X$. Does the oscillation
of $f$ decrease to zero on some sequence $E_n$
of infinite-dimensional subspaces of $X$?}

To state it in a more precise way we need some notation.
For a fixed function $f$ as above, and for an
arbitrary subspace $E \subset X$, let
$I_E(f) = [a(E), b(E)]$, where
$a (E)= \inf \{f(x) \mid x \in S \cap E\}$
and
$b (E)= \sup \{f(x) \mid x \in S \cap E\}$.
Then  let
$$
{\cal O}(f) = {\cal O}_X(f)=
 \inf \{ b(E) - a(E) \mid E \subset X, \dim E = \infty\}.
$$
The question then becomes:  is ${\cal O}(f)=0$?

If the answer is ``yes'' then
\begin{description}
\item[{}]
{\em there exists a real number $s$  such that for every $\ep >0$
there is a  subspace $E$ with $\dim E = \infty$ such that
$|f(x) - s| <\ep$ for all $x \in S\cap E$.}
\end{description}
The collection of all the numbers $s$ is called the {\em spectrum of $f$}
and denoted by $\gamma_\infty (f)$ or $\gamma_\infty (f, X)$ (see [M.69]).

And so, we are asking whether the spectrum $\gamma_\infty (f)$ is non-empty
for all uniformly continuous functions $f$ on the  sphere of an arbitrary
Banach space $X$, or  of some Banach space $X$?

Intuition says that the answer is {\em obviously} negative, at least for
$X = l_2$, say, because there is no reason for it to be positive. Uniform
continuity is a local geometric condition
with no  connection to a linear structure of a space,
and the existence of
$s \in \gamma_\infty (f)$ is a  global  linear property.
One never studies
what seems to be  obvious and the question was not 
an exception to this rule.

Note that James [J.64] showed that, in the above terminology,
$\gamma_\infty (f, l_1)$  and $\gamma_\infty (f, c_0)$ are non-empty
for $f$ being an equivalent norm on these spaces.
This result  did not
contradict the intuition, because the norms  in $l_1$ and $c_0$
are in a sense extremal, and the proofs deeply depended on this
fact. So, at the time, it did not even raise a similar
question for, say, $l_2$.

\subsection{}
\label{f_spectrum}

However, it was observed in 1967 ([M.67], \cf\ also [M.69], [M.71a])
that a slightly different {\em finite-dimensional spectrum}
$\gamma(f)$ is always non-empty. We say that $s \in \gamma(f)$
whenever
\begin{description}
\item[{}]
{\em for every $\ep >0$ and for every $n$ there exists an $n$-dimensional
subspace $E_n \subset X$ such that
$|f(x) - s| <\ep$ for all $x \in S\cap E_n$.}
\end{description}
We have the following fact valid
for every infinite-dimensional Banach space $X$.
\begin{fact}
  For  every uniformly
  continuous real function  $f$ on the unit sphere $S$,
  $\gamma(f)\ne \emptyset$.
\end{fact}

As we explained above, this somewhat contradicted
intuitions of that time.
Just to support  these  intuitions,
let us recall the Grinblatt's paper  [G.76]
where an example was presented
of a  bounded continuous, but not uniformly
continuous,  function $f$ on the sphere  $S$
in the Hilbert space,
which has the oscillation at least 1 on every 2-dimensional
central section of $S$.
So the  fact above indeed fundamentally rests on an interplay of uniform
continuity of a function and non-compactness of the sphere.

Thus, since the finite-dimensional spectrum $\gamma(f)$
involves subspaces of arbitrarily high dimensions
and it is always non-empty, it eventually  became
natural to expect that the (infinite-dimensional) spectrum
$\gamma_\infty(f)$ is also   non-empty.

\subsection{}
\label{tilda_sp}
Let us now consider the case when the function $f= \Snorm{\cdot}$
is another norm on $X$, continuous with respect to the original norm.
We have the following two mutually exclusive
possibilities.
\begin{description}
\item[(a) Spectrum:]
For every norm $f$ we have
$\gamma_\infty(f) \ne \emptyset$.
This would mean that either on some infinite-dimensional subspace
$f$ is arbitrarily small, if   $0 \in \gamma_\infty(f) $, or,
if $0\ne s \in \gamma_\infty(f) $, then $f$ is ``almost'' an isometry
on some infinite-dimensional subspace.
\item [(b) Distortion:]
There is a norm $f$
such that $\gamma_\infty(f) = \emptyset$.
This means that the norm $f$ has an oscillation
with respect to the original norm non-decreasing to zero
on {\em any}  infinite-dimensional subspace.
\end{description}

In view of the Fact above,
the existence of a norm satisfying  condition (b)
would be clearly connected with some very essential
infinite-dimensional  effects.

If a uniformly continuous function $f$ satisfies (b)
then, obviously, there exists an interval
$I = [\bt, \de]$, with $\bt < \de$, such that
\begin{description}
\item[(i)]
for every $ \ep >0$ there exists a subspace
$ Y= Y_\ep $, $\dim Y = \infty $,
such that
$I_Y(f) \subset [\bt-\ep, \de+\ep]$;
\item[(ii)]
For every
$ E \subset Y$, $\dim E = \infty $
one has
$ I_E(f) \supset (\bt, \de)  $.
(Here $Y$ is a subspace from (i) corresponding to
 $\ep = 1$, say.)
\end{description}
The collection of all such intervals $I$
is called the {\em tilda-spectrum} of $f$
and denoted by $\tilde{\gamma}(f)$.

Of course, the case $\de = \bt$ reduces the interval
to one point, $\bt \in \gamma_\infty (f)$, which we also
consider as a part of $\tilde{\gamma}(f)$.

Therefore we have (see [M.69])
\begin{fact}
    For  every uniformly
  continuous real function  $f$ on the unit sphere $S$,
  $\tilde{\gamma}(f)\ne \emptyset$.
\end{fact}

\subsection{}
\label{level_dist}
Note that  if  $f$ is a norm
on $X$ as in~\ref{tilda_sp}, and
if $I = [0, \de] \in \tilde{\gamma}(f)$,
then necessarily $\de =0$ (see [M.69]).

In the case $\bt >0$ we introduce a {\em level of distortion}
of an interval  $I \in \tilde{\gamma}(f)$
by $\dis (I)= \de / \bt$, and
a {\em level of distortion} of an equivalent norm
$f$ by
$$
\dis(f) = \sup \{ \dis(I)\mid I \in \tilde{\gamma}(f)\}.
$$

We have a similar  alternative as in \ref{tilda_sp}.
\begin{description}
\item[(a')] Either for any equivalent norm $f$ on $X$ one has $\dis (f)=1$,
\item[(b')] or there exists an
 equivalent norm $f$ on $X$ such that  $\dis(f)>1$.
\end{description}

In terms of the spectrum, condition (a') means that
for any equivalent norm $f$ on $X$ and  any
infinite-dimensional   subspace $Z$ of $X$,
the spectrum  $\gamma_\infty (f_{|Z})$
of the restriction of $f$  to $Z$, in non-empty.
Similarly, condition (b') means
that $X$ contains a distortable
infinite-dimensional subspace:
there  exists   an infinite-dimensional   subspace $Z$ of $X$
and an  equivalent norm $f$ on $X$
such that $\gamma_\infty (f_{|Z})= \emptyset$,
that is, $f$ is a distortion on $Z$.

It was proved by Milman in 1969   that
\begin{thm}
  Let $X$ be a Banach space. Assume that  $\dis (f) =1$
  for every equivalent
  norm $f$ on $X$. Then  either for some
  $ 1 \le p <\infty$, $X$ contains a   $(1+ \ep)$-isomorphic copy
  of $l_p$ (for every $\ep >0$), or  $X$ contains a
  $(1+ \ep)$-isomorphic copy of $c_0$ (for every $\ep >0$).
\end{thm}

(The result was  stated in [M.69], Section 3.3,  with  the
complete  proof   in [M.71b].)

And so, alternative (a') would imply an exciting structural
theory for Banach spaces.
However, in 1974, Tsirelson [Ts.74] constructed a space
$T$ which does not contain an isomorphic copy  of any $l_p$
($ 1 \le p < \infty$) or of $c_0$. This means that the space
$T$ satisfies the alternative (b'): $T$ contains a distortable
infinite-dimensional  subspace $Z$.
(In fact, it can be shown by a direct argument
that $T$ itself is also distortable.)

An interesting feature of Tsirelson's example is that the norm
is not given by an explicit formula but it is defined by
an equation. This was the first construction of such a type,
and  essentially, with only minor modifications, the only one.
In the dual form, which has been put forward by Figiel and
Johnson [F-J.74], the norm is defined, for
a finite sequence of real numbers
$x \in \RRN$, by
\begin{equation}
\|x\|_T = \max \, \left\{ \|x\|_{c_0}, {\frac{1}{2}}
   \sup
 \sum_{i=1}^n \|E_i x\|_T \right\},
  \label{tsir}
\end{equation}
where the inside supremum is taken over all succesive
intervals $\{E_i\}$ of positive integers such that
$n < \min E_1 \le \max E_1  < \min E_2 \le \dots
< \max E_{n-1}< \min E_n$
and over all $n$. For $x = \sum_i t_i e_i \in X$ and
an interval $E$, we set  $Ex = \sum_{i\in E} t_i e_i $.
Tsirelson's space $T$ is then a completion of
$ \RRN$ under the norm $\|\cdot\|_T$.
Most of important properties of the space $T$ and
related spaces can be found in [C-S.89]
and references therein.

\subsection{}
\label{as_s}
Let us return to a distortion situation
when  $I = [\bt, \de] \in \tilde{\gamma}(f)$
with  $\bt < \de $ and let us give its geometric
interpretation.

Let $\ep < (\de - \bt)/2$ and let $Y = Y_\ep$ be a corresponding
subspace. Define two sets
$$
A = \{ x \in S\cap Y \mid f(x) < \bt + \ep\}
\quad \mbox{\rm and}\quad
B = \{ x \in S\cap Y \mid f(x) > \de - \ep\}.
$$
For every infinite-dimensional subspace $E \subset Y$
we have $A \cap E \ne \emptyset$
and $B \cap E \ne \emptyset$.
A set satisfying such a property is
called an \as\ set (in $Y$).
So in our situation, $A$ and $B$ are two
\asm\ sets with positive distance
apart, $\dist(A, B)>0$.
The fact of
the existence of such a pair $(A, B)$
is thus a consequence of distortion.
Conversely, this fact also implies some distortion
property. The Urysohn function
for sets $A$ and $B$ is a uniformly continuous function
with an empty spectrum $\gamma_\infty(f)$; to construct
an equivalent norm without spectrum some
additional convexity assumptions are required.

\subsection{}
\label{o_s}
Given a distortion situation it is natural to ask
a quantitative question, how large can  $\dis (f)$ be.
Note that Theorem~\ref{level_dist} ensures only the existence
of distortion  but provides no quantitative information
on  $\dis (f)$.

It is to Rosenthal's credit that in 1988 he
asked the first named author, Odell and several others,
how to find a direct formula for a distortion
on Tsirelson's space, and how large such a distortion can be.
Odell (unpublished) in 1989/90 constructed two \asm\
sets in $T$. He also showed that the spaces $T_{\lambda}$,
obtained by replacing $1/2$ in
the definition (\ref{tsir}) by $1/\lambda$,
have distortions  $d_\lambda  $
of order $1/\lambda$, hence
$d_\lambda  \to \infty$
as $\lambda \to \infty$.
So, for every real number $d$, there is
a space with  a level of distortion at least $d$.

Let us mention that  an approach to distortions
using the theory of Krivine--Maurey types
was presented in [H-O-R-S.91]. In particular this paper
contains another proof of  Theorem~\ref{level_dist}.

The next step was done by Schlumprecht [S.91],
who changed $1/2$ to $1/ \ln n$,
which also allowed to
start $E_1$ at any place (not necessarily far out).
This had an important
effect on the geometry of the space:
the unit vector basis becames subsymmetric
and, as Schlumprecht showed, any distortion level
is attained by some equivalent norm.
Schlumprecht's space $S$ also has the
property of an {\em infinite distortion}:
\begin{description}
\item[{\ \  }]
{\em there exists a sequence of asymptotic
sets $\{A_i\}$ on the  sphere of $S$ \newline
such that $\dist(A_i, \conv(\bigcup_{j\ne i} A_j))
\ge 1$.}
\end{description}
In fact, $S$ satisfies still stronger condition that
there also exists  a sequence of
sets $\{A_i^*\}$ on the  sphere of the dual space $S^*$
such that the system $\{A_i, A_i^*\}$ is ``nearly
biorthogonal''.

This was the starting point for Gowers' and Maurey's
construction.

Finally, this year, Odell and Schlumprecht [O-S.92]
proved that for every $1 < p < \infty$, $l_p$ has
an arbitrarily large (and even infinite) distortion, this way
finishing off the problem which originated from [M.69],
[M.71b].
Again, they did not construct
asymptotic sets far apart in, say, $l_2$, but transformed them
in an ingeneous  non-linear way from Tsirelson's space, or, on more
advanced level, from Schlumprecht's space.
Combining this outstanding result with Theorem~\ref{level_dist}
we see that

\begin{thm}
  Any Banach space $X$ which does not  hereditarily contain copies
  of $l_1$ and $c_0$, contains a  distortable subspace,
  \ie\ there exists an equivalent norm $f$
  on  $X$ such that $\dis (f) >1$.
\end{thm}

Moreover, Odell and Schlumprecht proved that
on the sphere  $S(l_1)$  there is
a Lipschitz  function $f$ (not a norm)
with an empty spectrum, $\gamma_\infty (f, l_1) = \emptyset$.
It was shown earlier by Gowers [G.91] that
 $\gamma_\infty (f, c_0) \ne \emptyset$, for any uniformly
continuous function on $S(c_0)$.

\subsection{}
\label{general_dist}
Let us go back to the quantitative  question: does any Banach space
$X$ not containing  hereditarily
copies of $l_1$ and $c_0$, have   arbitrarily large  distortions?

This is not yet clear. To study this problem,
we consider in this paper spaces with
{\em bounded distortions}.
These are spaces $X$ such that
for some constant $D$ we have
$\dis(f) \le D$, for every infinite-dimensional
subspace $Z$ of $X$ and every equivalent norm
$f$ on $Z$.
What kind of simple ``basic structural blocks''
(\ie\ subspaces) can
such a space $X$ contain?
To explain our result let us define the class
of \as\ $l_p$ spaces.
(The rather standard notation concerning successive blocks
of a basis and related concepts
will be explained at the beginning of the next section.)

\medskip
\noindent{\bf{Definition\ \ }}
A Banach space $X$ with  a normalized basis $\{x_i\}$
is said to be \as\ $l_p$ space, for some
$1 \le p < \infty$ (resp. \as\ $c_0$ space)
if  there exists a  constant $C$ such that
for every $n$ there exists
$N= N(n)$ such that any normalized successive blocks
$N <z_1 < z_2 < \ldots< z_n$ of $\{x_i\}$ are $C$-equivalent
to the unit vector basis in $l_p^n$ (resp.\ in $l_\infty^n$).
By ${\la}_p (X)$ we denote the infimum of
all constants $C$ as  above.

\medskip

Note that Tsirelson's space $T$ is an \as\ $l_1$ space which does not
contain a subspace  isomorphic to $l_1$.

For spaces with bounded distortions, let
\begin{equation}
\dis (X) = \sup  \dis (f),
  \label{dis}
\end{equation}
where the supremum is taken over all
equivalent norms $f$ on $X$.

Recall a standard and easy observation that if
$Z \subset X$ is an infinite-dimensional subspace
and $f$ is an equivalent norm on $Z$ then
there exists an equivalent norm $\tilde{f}$ on $X$
such that $\tilde{f}_{|Z} = f$. This immediately
implies that $\dis (Z) \le \dis (X)$.

\begin{thm}
  Let $X$ be a Banach space  with   bounded distortions
  and let $\dis(X) < D$.
  There exists a
  subspace $Y$ of $X$ which is either \asm\ $l_p$, for some
  $1 \le p < \infty$, or  \asm\ $c_0$. Moreover,
  ${\la}_p (Y)$ depends on $D$ only.
\end{thm}

We learned recently that B.~Maurey [Ma.92] also proved
this theorem and used it to show that every space
of  type $p >1$ with an unconditional basis has
arbitrarily large distortions.

In contrast with the result for $  \dis (X) =1$
(Theorem~\ref{level_dist}),
the theorem above  recognizes, as ``basic structural blocks'',
a {\em class} of Banach spaces rather than a concrete space,
as it was suggested by a ``naive'' intuition of the 60s.
(It is well-known that
varying $\lambda$ in the definition of Tsirelson's spaces
$T_\lambda$ we get a sequence of non-isomorphic \asm\ $l_1$
spaces, and the so-called $p$-convexified  Tsirelson's spaces
show that the same phenomenon holds
for any fixed $1 \le p < \infty$ or $c_0$.)

Another important point is a difference with the local theory
of Banach spaces. The definition of \asm\ $l_p$ spaces
is ``almost''
local, in that it involves finite-dimensional subspaces
and parameters not depending on the dimension.
However, this isomorphic
definition does not imply a  $(1+\ep)$-isometric version,
in the standard spirit of the local theory.
Indeed, we would call a Banach space
an \ai\ \as\ $l_p$ space if
$ \la_p (X)=1$.  It is well-known to specialists
that every
  \ai\ \as\ $l_p$ space contains,
for every $\ep >0$, a subspace  $(1+\ep)$-isomorphic
to $l_p$ (see~\ref{isomorph} for a short argument).

Our method involves geometry of infinite-dimensional sphere
and  a suitable geometric language will be introduced
in the next section.

\section{Preliminaries}
\subsection{}
\label{basis}
Since in this paper we are concerned with the existence
of nice infinite-dimensional subspaces
inside Banach spaces from a certain class,
we may and will assume, unless  stated otherwise,
that Banach spaces discussed here have a monotone
basis.
In such a situation we will use the standard notions
of the dual basis, equivalent bases, block bases,
block subspaces,  basic sequences, etc.
They can be found \eg\ in [L-T.77].
Let us only mention that we will say that two basic sequences
$\{x_i\}$ and $\{e_i\}$ are  $C$-equivalent,
for some constant $C$,
if for any (finite) sequence of scalars $\{a_i\}$
we have
$$
C^{-1} \bigl\|\sum_i a_i x_i \bigr\| \le \bigl\|\sum_i a_i e_i \bigr\|
\le C \bigl\|\sum_i a_i x_i \bigr\|.
$$
We will consider only vectors with finite
support.
For vectors $x, y \in X$ and  subspaces $E, E_1, E_2$,
we will freely use the
notation  $ n < x$ to denote that
$n < \min \supp (x)$; then
$x < y $ if  $\max \supp(x) < \min \supp (y)$;
then  $x < E$ if
$x < y$ for every $y \in E$;
and  $E_1 < E_2$ if
$x < y$ for every $x \in E_1$ and $ y \in E_2$.

For a Banach space $X$
by $B_X$ and  $S(X)$  we denote the unit
ball and the unit sphere in $X$, respectively;
for a subspace $E \subset X$ we set
$B_E= B_X \cap E$ and $S(E)= S(X) \cap E$.

\subsection{}
\label{as_s2}

Asymptotic sets were defined in~\ref{as_s}
where their  basic connection to distortions
was indicated.
To get a better  understanding of their geometric  properties
let us make some  easy general  observations,
valid for arbitrary Banach spaces
(which may have no basis).

\begin{fact}
  Let $(Z, \|\cdot\|)$ be a Banach space and let $\Snorm{\cdot}$
  be a seminorm on $Z$ such that  $\Snorm{z} \le  \|z\|$,
  for all $z \in Z$.  Assume that there exists
  an asymptotic set $A \subset S(Z, \|\cdot\|)$ such that
  $\|\cdot\|$ and  $\Snorm{\cdot}$ are equivalent on $A$.
  Then there exists a subspace $E $ of $Z$ of finite-codimension
  such that $\Snorm{\cdot}$ is a norm on $E$ equivalent to
  $\|\cdot\|$.
\end{fact}
\proof
Clearly,   $\Snorm{\cdot}$ is a norm on the subspace $W$
spanned by the set $A$, and since $A$ is \asm\ then
$\codim W <\infty$. A standard well-known fact
(\cf\ \eg\ [K.66], [L-T.77]) implies that if
$\Snorm{\cdot}$ and $\|\cdot\|$ were not equivalent
on any subspace $E $ of $W$ of finite-codimension
then for every $\ep >0$ there would be an infinite-dimensional
subspace $F$ of $W$ such that  $\Snorm{z} \le \ep  \|z\|$,
for all $z \in F$. But for $\ep$ sufficiently
small this is impossible, since  $F$ intersects $A$.
\qed

\medskip
\rem
Let   $Z$ be a Banach space with   bounded distortions,
  $\dis(Z)< D$,
and  let $A \subset S(Z)$ be  an \asm\ set symmetric about the origin.
Then there exists an infinite-dimensional subspace
$ F$ of $Z$ such that
$(1/D) (B_Z \cap F) \subset \conv A$.

Indeed,
let  $W$ be the finite-codimensional
subspace spanned by $A$.
Applying the fact above to the norm  $|\cdot|$ on $W$
whose   unit ball is  $\conv A \cap W$,
we get that  $\|\cdot\|$
and $|\cdot|$ are equivalent on a certain
subspace $E$ of $W$ of finite codimension.
Thus   there exists an infinite-dimensional
subspace $F$ of $E$
such that
\begin{equation}
\|z\| \le |{z}| \le D \|z\|
\qquad \mbox{\rm for} \qquad  z \in F.
  \label{as_xxx}
\end{equation}
This implies the required inclusion
of the corresponding unit balls.

\subsection{}
\label{BBB}
To make the arguments more compact,
we introduce several short
notations for certain families
of subspaces of a given Banach space $Z$
(with a basis).
Typically, $Z$ will be a block subspace
of the fixed Banach space $X$.

By  ${\cB}_\infty(Z)$ we
denote the  family of all
infinite-dimensional block subspaces $E \subset Z$;
next,  ${\cB}^t(Z)$  denotes the  family of all
(block) subspaces $E \in \cB_\infty (Z)$
of finite-codimension, \ie\
$\dim Z/E <\infty$; finally, if
$Y \in \cB_\infty (Z)$ and $z \in S(Z)$
(with finite support), then  ${\cB}^t(Y,z)$
denotes the family of all subspaces $F \in \cB_\infty (Y)$
such that $z < w$ for all $w \in F$.

\subsection{}
\label{moduli}
Let us recall the geometric notions  of  asymptotic
averages and moduli, which play a major role in
our approach. These notions
were introduced and studied by Milman in 1967--70.
A survey on this subject can be found
in [M.71b], \cf\  also more recent paper [M-P.89].

The moduli are defined relatively to a fixed
family $\cB$ of subspaces of a space $X$,
which satisfies the
filtration condition
 $$
\mbox{\rm For\ every\ } E_1, E_2 \in {\cB}
  \mbox{\rm \  there\  exists \ }
  E_3 \in {\cB} \mbox{\rm \  such\  that\ } E_3 \subset E_1 \cap E_2.
 $$
Typically, the family $\cB$ will be
$\cB^t$, which have been  defined in~\ref{BBB},
although we will make an exception from this rule
in Section~\ref{non-dist}.

For a continuous bounded function $h: S(X) \to \Rn{}$
and an infinite-dimensional  subspace $E\subset X$
define lower and upper moduli $\bt$- and $\de$-,
respectively, by
\begin{eqnarray}
\bt [h, \cB, E] &=& \bt_x [h, \cB(E)]=
 \sup_{F \in \cB(E)}\
  \mathop{\vphantom{p}\inf}_{x \in S(F)} h(x),\nonumber\\
\de [h, \cB, E] &=& \de_x [h, \cB(E)] =
\mathop{\vphantom{p}\inf}_{F \in \cB(E)}\
    \sup_{x \in S(F)} h(x).
\label{loc_mod}
\end{eqnarray}
For a continuous bounded function $f: S(X)\times S(X) \to \Rn{}$
and $E \in \cB_\infty(X)$  we set
\begin{eqnarray}
\bt\bt [f, \cB, E]=
    \bt_x \left[\bt_y [f(x,\cdot), \cB(E,x)], \cB, E\right]\nonumber\\
\de\de [f, \cB, E]=
    \de_x \left[\de_y [f(x,\cdot), \cB(E,x)], \cB, E\right].
\label{glob_mod}
\end{eqnarray}

Observe that if $E \subset F$ then
$$
\bt [h, \cB, F] \le \bt [h, \cB, E] \le \de [h, \cB, E]
\le \de [h, \cB, F],
$$
hence also
$$
\bt\bt [f, \cB, F]\le \bt\bt [f, \cB, E]\le
\de\de [f, \cB, E] \le \de\de [f, \cB, F],
$$
for all functions $h$ and $f$ as above.

\section{Non-distortable spaces}
\label{non-dist}

To develop better geometric intuitions
and to illustrate the use of the
$\bt$- and $\de$- averages
we start with
the  isometric case and
we will sketch the proof  of
Milman's
theorem on non-distortable spaces,
Theorem~\ref{level_dist}.

\subsection{}
\label{non_notation}
In the isometric situation discussed here  there
is  no real advantage in passing to a subspace with a basis,
in fact, this would  confuse a geometric picture rather
than clarify it. Therefore we present an argument
which makes no reference to the existence of  a basis
and thus it works for  an arbitrary Banach space.
The averages we will consider here will be taken
with respect to the family
$\cB= \cB^0(E)$ of all finite-codimensional subspaces
of a given space $E$.
This family clearly satisfies
the filtration  condition.

We will consider
the collection of functions
\begin{equation}
      f_\ep (x,y) = \|x + \ep y\|-1
      \quad \mbox{\rm  for} \quad x, y \in S(X).
  \label{f_ep}
\end{equation}

The averages of functions $ f_\ep$ will be called
the $\bt$- and $\de$-moduli, as they reflect
a geometric behaviour of the sphere in a Banach
space.

For a subspace $E \subset X$ and $x \in S(X)$,
the notation $\cB(E,x)$ used in (\ref{glob_mod})
simply means $\cB^0(E)$.
Also, $\cB_\infty(E)$ denotes the family of all
infinite-dimensional subspaces of $E$.

We  will consider the local modulus
$\bt_y [f_\ep (x,\cdot), \cB^0(E,x)]$
denoting it by $\bt(\ep, x, E)$
and the global modulus
$\bt\bt [f_\ep, \cB^0, E]$, denoting it
by  $\bt \bt(\ep, E)$. Similarly, the local modulus
$\de_y [f_\ep (x,\cdot), \cB^0(E,x)]$
will be denoted by $\de(\ep, x, E)$
and the global modulus
$\de\de [f_\ep, \cB^0, E]$,  by $\de \de(\ep, E)$.

To illustrate the expected behaviour of the moduli,
let us observe that for $X = l_p$, $ 1 \le p < \infty$
we have
$$
\bt(\ep, x, l_p) = \de(\ep, x, l_p) = (1 + \ep ^p)^{1/p} -1,
$$
for all $x \in S(l_p)$ and all $\ep >0$. This function has the
order $\ep ^p /p$ as $\ep \to 0$.

Computation of the moduli for some other spaces can be found in
[M.71b].

\subsection{}
\label{bb=dd}
\begin{lemma}
  Let $X$ be a Banach space such that  $\dis (X) =1$.
There exists an infinite-dimensional subspace $F$ of $X$ such
that
  \begin{equation}
    \bt \bt(\ep, F)= \de \de(\ep, F)
    \quad \mbox{\rm for}\quad  \ep >0.
    \label{stab_0}
  \end{equation}
\end{lemma}
\proof It is not difficult to see that if $\dis (f) =1$ for
every equivalent norm $f$ on $Z$ then ${\cal O}_Z (g)=0$ for
every $Z \in \cB_\infty (X)$ and every uniformly continuous
convex function $g: Z \to \Rn{}$.  We will show this at the
end of the proof.

Observe that for each $\ep >0$ and $x \in S(X)$, the
function $f_\ep (x, \cdot)$ is convex, therefore
${\cal O}(f_\ep (x, \cdot))=0$.
Stabilizing over $y$
with a given $\vt >0$  we get a subspace
$\tilde{E} \in \cB_\infty(X)$ such that
\begin{equation}
  0 \le \sup_{y \in S(\tilde{E})} f_\ep (x,y) -
  \mathop{\vphantom{p}\inf}_{y \in S(\tilde{E})} f_\ep (x,y) <
  \vt \quad \mbox{\rm for} \quad x \in S(\tilde{E}).
  \label{stab_1}
\end{equation}

Now we take a dense set $\{x_i\}$ in the unit sphere $S(X)$
and a sequence of $\vt_i \downarrow 0$, and we let, for every
$i= 1, 2, \ldots$, $\ep$ vary over a finite $\vt_i$-net
${\cal N}_i$ in $[\vt_i, 1/ \vt_i]$; this way we can construct
a sequence $E_1 \supset \ldots \supset E_i \supset\ldots$ of
infinite-dimensional subspaces such that on $S(E_i)$ the
inequality analogous to (\ref{stab_1}) holds for $x_j$, with
$j = 1, \ldots, i$,
and $\vt_i$ and  all $\ep \in {\cal N}_i$, ($i=1, 2, \ldots$).
Picking $e_i \in S(E_i)$
for $i= 1, 2, \ldots$ and setting
$E = \spn [e_i]$, we get,
by this diagonal procedure,
$E \in \cB_\infty (X)$
such that for all $x \in S(E)$ and $\ep >0$
we have
$$
\bt(\ep, x, E)= \de(\ep, x, E).
$$

Now, with a fixed $\ep >0$, the function $\de(\ep, \cdot, E)$
is again convex. Stabilizing over $x$, with a fixed $\ep >0$,
and then passing to a diagonal in a similar
way as before,  we get an
infinite-dimensional  subspace $F$ of $E$
on which (\ref{stab_0}) holds.

It remains to show that
${\cal O}_Z (g)=0$ for every $Z \in \cB_\infty (X)$ and every
uniformly continuous convex  function $g: Z \to \Rn{}$.
To this end, fix  $Z \in \cB_\infty (X)$ and
pick $[\bt, \de]\in \tilde{\gamma}(g)$,
and,  for an  arbitrary (fixed) $\ep >0$,
let $Y \in \cB_\infty (Z) $ be a corresponding stabilizing subspace,
satisfying conditions (i) and (ii) from \ref{tilda_sp}.

Assume first that $g(x) = g(-x)$ for $g \in S(X)$.
Consider the symmetric
asymptotic sets $A$ and $B$
defined in~\ref{as_s}. By Remark in~\ref{as_s2} (with $D =1+\ep$)
we get  a subspace
$E \in \cB_\infty (Y)$ such that
$B\cap E \subset B_E = B_X \cap E  \subset (1+\ep)\conv A$. Thus,
by  convexity and uniform continuity
of $g$, we have
$$
 \de - \ep - \vt \le \inf_{z \in B\cap E} g(z)
\le \sup_{w \in A} g(w)+ \vt \le \bt + \ep + \vt,
$$
where
$  \vt = \vt(\ep)\to 0$  as $\ep \to 0$.
Letting $\ep \to 0$ we get $\bt = \de$,
hence ${\cal O}_Z(g) =0$.

If $g$ is arbitrary,
set $h (x) = (1/2)(g(x)+g(-x))$.
Since  ${\cal O}_Y (h) =0$,  find
a subspace
$Y_1  \in \cB_\infty (Y)$ such that
$|h(x) - s| < \ep$
for all $x \in S(Y_1)$.
If  $\ep $ is small enough
then
$|s -  (\bt + \de)/2 | < 2\ep$,
hence $\bt + \ep < s < \de - \ep$.
(In fact, we will use the later inequality
only.)
Indeed, otherwise $g(x)$ and $g(-x)$
would not  compensate each other.
Formally,
pick $w, v \in S(Y_1)$ such that $g(w)= \bt + \ep$ and
$g(v)= \de - \ep$ and observe that
$$
s -(\bt + \de)/2 \le \bigl(s - h(w)\bigr) +
\bigl( (g(w)+g(-w))/2 -   (\bt + \de)/2 \bigr) \le 2 \ep;
$$
and similarly, using $h(v)$, it is easy to establish
the lower estimate by
$- 2 \ep$.

Consider the set
$ A = \{y \in S(Y_1) \mid |g (y)- s| \le \ep \} $,
which is asymptotic in $Y_1$.
Observe that for $y \in A$ we have
$$
|g (-y)- s| \le |g (y)+ g (-y)- 2s| + |g (y)- s| \le 3\ep,
$$
so that $A$ is ``almost'' symmetric.
By Remark in \ref{as_s2}
we get  a subspace
$E  \in \cB_\infty (Y_1)$ such that
$ B_E \subset (1+\ep)\conv (A \cup -A)$. Thus
for every $z \in S(E)$ we have,
$$
 g(z) \le \sup_{w \in A \cup -A} g(w)+ \vt \le s + 3\ep + \vt,
$$
where
$  \vt = \vt(\ep)\to 0$  as $\ep \to 0$.
This, combined with the
stabilization property of $h$,
yields that ${\cal O}_Z(g)$ is arbitrarily small,
hence equal to 0, as required.
\qed

\medskip
\rem
In the non-distortion situation
of the lemma it is not difficult to show
that 
$X$ contains a  subspace with a basis
$\{u_i\}$ such that for every $n$ and
for all blocks $ n < w < v$
we have $\max\,( \|w + v\|, \|w - v\| )
\le (1+ 2^{-n}) \|w + v\|$.
A standard argument shows that the tails of $\{u_i\}$
are unconditional with the constants as close to 1 as we wish.
We could then   consider
the moduli related  to the family $\cB = \cB^t$
and  note that the argument from
\ref{sketch} applies for these moduli as well.
Moreover, in this  situation it would be clearly sufficient
to discuss the observation  opening the  proof
of the lemma only  for symmetric functions.

\subsection{}
\label{sketch}
{\bf Sketch of the proof of Theorem~\ref{level_dist}}\ \
We will show that the equality
$ \bt \bt(\ep, F)= \de \de(\ep, F)$ yields that
the space $F$  contains
$(1+\ep)$-isomorphic copies of $l_p$ or $c_0$.

Using definitions of
$ \bt \bt(\ep, F)$ and $ \de \de(\ep, F)$ it is
possible, given $\eta > 0$, to construct a basic
sequence $\{y_i\}$ in $F$ such that for
any finite sequence  $\{a_i\} \in \RRN$ of real numbers
the norm
$\| \sum_i a_i y_i \|$
admits an upper estimate by
$(1+\eta) \Phi(\{\de \de(|a_i|,F)\})$
and a lower estimate by
$(1+\eta)^{-1} \Phi(\{\bt \bt(|a_i|,F)\})$,
where $\Phi(\{\cdot\})$ is a real function
defined on the space of all finite sequences
of real numbers (see [M.71b], Theorem 4.5).
Moreover, the same estimates are satisfied
for every vector of the form
$ \sum_i a_i u_i $,
where  $\{u_i\}$ is  a block
basis of  $\{y_i\}$.
(This  argument is similar to
a well-known construction sketched in~\ref{isomorph}.)
By Lemma~\ref{bb=dd}, this implies
that all block bases of
$\{y_i\}$ are $(1+\eta)$-equivalent.
By Zippin's theorem, the
basis   $\{y_i\}$ is
$(1+\eta)^{\alpha}$-equivalent
to the unit vector basis in $l_p$,
for some $1 \le p <\infty$ or in  $c_0$
(here $\alpha>0$ is a numerical constant).
\qed



\section{Tilda-spectrum in general}
\subsection{}
\label{intuition}
Let $X$ be a Banach space with a basis $\{x_i\}$.
Let $f(z_1, z_2, \ldots, z_l)$ be a
uniformly continuous real function defined
on  sequences of $l$ normalized block vectors
$z_1 < z_2< \ldots < z_l$.
First let us describe a rough intuition 
of an interval $[\bt, \de]$ of tilda-spectrum
of $f$ on a subspace $\tilde{Y} \in \cB_{\infty}(X)$,
leaving the precise definition
for  later parts of this section.

Let $\tilde{I}$ be  the closure
of the interval
of values of $f$ on $\tilde{Y}$.
By restricting the domain of the variable
$z_l$  to any
subspace of $\tilde{Y}$,
with other variables fixed,
we do not increase $\tilde{I}$.
Therefore, for any $z_1 < z_2< \ldots < z_{l-1}$
fixed,
let  $Y_1$ be a subspace of $\tilde{Y}$
such that passing with $z_l$ to  $Y_1$
corresponds to the ``maximal decrease''
of $\tilde{I}$.
Let $I^{(1)}$ be the closure of the interval of values
of this restricted  $f$.
Continue
the procedure of restricting $z_{l-1}$,
with  $z_1 < z_2< \ldots < z_{l-2}$ fixed.
The closed interval $ I^{(l)}= [\bt, \de]$
obtained after the $l$-th step
is called an interval
of the tilda-spectrum of $f$
in $\tilde{Y}$.

\subsection{}
\label{def}
Let  $\tilde{Y} \in \cB_{\infty}(X)$.
The precise definition of the tilda-spectrum
of $f$ on  $\tilde{Y} $
involves the notions of the
$\bt$- and $\de$- averages, introduced
in~\ref{moduli}.
These averages will be applied to
functions of the form
$h (z_1, \ldots, z_k)$, considered as
functions of  $z_k$
with $z_1< \ldots< z_{k-1}$ fixed,
and with respect to
the   family
$\cB^t (Y,z_{k-1})$ of all finite-codimensional
block subspaces of $Y$ with the support after $z_{k-1}$
($Y \in \cB_\infty (X)$ is a subspace).
To make the formulas more compact,
we will  indicate the variable $z_{k}$ and the subspace $Y$
in the subscripts, leaving  $z_{k-1}$
to be understood from the context. Thus we will
write, {\it e.g.\/} $\bt_{z_k,Y}(h(z_1,  \ldots,z_{k-1}, z_k) )$
for  $\bt[h(z_1,  \ldots,z_{k-1},\cdot), \cB^t (Y,z_{k-1})]$,
and so on.

We say that an interval
$[\bt, \de]$
is in the tilda spectrum of $f$
on  $\tilde{Y} $
if  there is  a subspace
$Y \in \cB_\infty(\tilde{Y})$
such that
the following two conditions are satisfied:
\begin{description}
\item[(i)]
$\bt = \bt_{z_1,Y}\left( \bt_{z_2,Y}
      (\ldots (\bt_{z_l,Y}(f(z_1, z_2, \ldots, z_l) ))
    \ldots)\right)$\newline
and \\
$\de = \de_{z_1,Y}\left( \de_{z_2,Y}
      (\ldots (\de_{z_l,Y}(f(z_1, z_2, \ldots, z_l) ))
     \ldots)\right)$;
\item[(ii)]
for all $H_1, H_2, \ldots, H_l \in  \cB_ \infty(Y ) $,
each of the averages $ \bt_{z_i,Y} $ and $ \de_{z_i,Y} $
in (i) can be replaced by
$ \bt_{z_i,H_i} $ and $ \de_{z_i,H_i} $, respectively;
that is, we have\\
$\bt = \bt_{z_1,H_1}\left( \bt_{z_2,H_2}
      (\ldots (\bt_{z_l,H_l}(f(z_1, z_2, \ldots, z_l) ))
    \ldots)\right)$\newline
and \\
$\de = \de_{z_1,H_1}\left( \de_{z_2,H_2}
      (\ldots (\de_{z_l,H_l}(f(z_1, z_2, \ldots, z_l) )  )
      \ldots)\right)$.
\end{description}

A subspace $Y$
for which the above conditions
hold is
called a spectrum subspace corresponding
to $[\bt, \de]$.  Then
any further subspace $Y'$ of $Y$
is a spectrum subspace as well.

\subsection{}
\label{new_mod}
To prove the existence of the tilda-spectrum
defined in~\ref{def}, it is convenient to introduce
modified averages $\bt^{st}$ and  $\de^{st}$.
The definition requires several steps.

\subsubsection{}
\label{new_mod_1}
Let $h (z_1, \ldots, z_k)$ be a uniformly continuous
function and let $E \in \cB_\infty (X)$.
Fix normalized blocks $z_1<  \ldots < z_{k-1}$
and set
$$
\gamma(E) = \inf  \bigl(
\de_{z_k,G} (h(z_1,  \ldots,z_{k-1}, z_k)) -
\bt_{z_k,G} (h(z_1,  \ldots,z_{k-1}, z_k)) \bigr),
$$
where the infimum is taken over all subspaces
$G \in \cB_\infty (E)$.
Pick $\ep_i \downarrow 0$ and construct a sequence
$G_0 = E \supset G_1 \supset \ldots$, such that
$G_i \in \cB_\infty (G_{i-1})$
and
\begin{equation}
\de_{z_k,G_i} (h(z_1,  \ldots,z_{k-1}, z_k)) -
\bt_{z_k,G_i} (h(z_1,  \ldots,z_{k-1}, z_k))
\le \gamma(G_{i-1}) +\ep_i,
\label{111}
\end{equation}
for $i = 1, 2, \ldots$.
Set
\begin{eqnarray}
\bt_{z_k,E}^{st} (h) &=&
\bt_{z_k,E}^{st} (h; z_1,  \ldots,z_{k-1}))\nonumber\\
&& \qquad = \lim_{i \to \infty}
\bt_{z_k,G_i} (h (z_1,  \ldots,z_{k-1}, z_k)),\nonumber \\
\de_{z_k,E}^{st} (h) &=&
\de_{z_k,E}^{st} (h; z_1,  \ldots, z_{k-1})) \nonumber\\
&& \qquad =  \lim_{i \to \infty}
\de_{z_k,G_i} (h (z_1,  \ldots,z_{k-1}, z_k)).
  \label{112}
\end{eqnarray}

Let $G = \spn [u_i]$ be a diagonal subspace for $\{G_i\}$,
that is,
$u_i \in G_i$ for $i = 1, 2, \ldots$. It is easy to
check from (\ref{111}) that for every
subspace $H \in \cB_\infty (G)$ we have
\begin{equation}
\bt_{z_k,E}^{st} (h) = \bt_{z_k,H} (h)
\qquad \mbox{\rm and} \qquad
\de_{z_k,E}^{st} (h) = \de_{z_k,H} (h).
  \label{113}
\end{equation}

\subsubsection{}
\label{new_mod_2}
Now let $\{(z_1^{(i)}, \ldots, z_{k-1}^{(i)}) \}$
be a dense countable subset of $(k-1)$-tuples
of  normalized blocks $z_1<  \ldots < z_{k-1}$.
Let $G^{(1)}=G$ be the subspace constructed at the end
of~\ref{new_mod_1} for $(z_1^{(1)}, \ldots, z_{k-1}^{(1)})$.
Starting from $G^{(1)}$,
construct $G^{(2)} \in \cB_\infty (G^{(1)})$
for the tuple $(z_1^{(2)}, \ldots, z_{k-1}^{(2)})$.
Proceeding by induction and using (\ref{113})
we get a sequence
of subspaces
$E \supset G^{(1)} \supset G^{(2)}  \supset\ldots$
such that for every $i= 1, 2, \ldots$
and all $H, H' \in \cB_\infty (G^{(i)})$
we have
\begin{eqnarray*}
\bt_{z_k,H} (h (z_1^{(i)},  \ldots,z_{k-1}^{(i)}, z_{k}))
&=& \bt_{z_k,E}^{st} (h; z_1^{(i)},  \ldots,z_{k-1}^{(i)})\\
&&\qquad = \bt_{z_k,H'} (h (z_1^{(i)},  \ldots,z_{k-1}^{(i)}, z_{k})),\\
\de_{z_k,H} (h (z_1^{(i)},  \ldots,z_{k-1}^{(i)}, z_{k}))
&=& \de_{z_k,E}^{st} (h; z_1^{(i)},  \ldots,z_{k-1}^{(i)})\\
&&\qquad = \de_{z_k,H'} (h (z_1^{(i)},  \ldots,z_{k-1}^{(i)}, z_{k})).
\end{eqnarray*}

Taking once more a diagonal subspace we get
$F = \spn [v_i]$, with $v_i \in G^{(i)}$ for $i = 1, 2, \ldots$
such that for all $(k-1)$-tuples
$(z_1, \ldots, z_{k-1})$
and for all subspaces
$H \in \cB_\infty (F)$ we have
\begin{eqnarray}
\bt_{z_k,H} (h (z_1,  \ldots,z_{k-1}, z_{k}))
&=& \bt_{z_k,E}^{st} (h; z_1,  \ldots,z_{k-1})\nonumber\\
\de_{z_k,H} (h (z_1,  \ldots,z_{k-1}, z_{k}))
&=& \de_{z_k,E}^{st} (h; z_1,  \ldots,z_{k-1}).
  \label{114}
\end{eqnarray}

\subsubsection{}
\label{new_mod_3}
Coming back to the definition~\ref{def} of the
tilda-spectrum, fix
a function $f = f(z_1, \ldots, z_l)$
and   a subspace  $\tilde{Y} \in \cB_\infty (X)$.

The existence of a tilda-spectrum interval
$[\bt, \de]$ will be proved by providing
explicit formulae for $\bt$ and $\de$
in terms of the stabilized
averages $\bt^{st}$ and $\de^{st}$.
This is done by the backward induction.

Let
\begin{eqnarray*}
  b_1 (z_1,  \ldots,z_{l-1})
&=& \bt_{z_l,\tilde{Y}}^{st} (f; z_1,  \ldots,z_{l-1}) \\
  d_1 (z_1,  \ldots,z_{l-1})
&=& \de_{z_l,\tilde{Y}}^{st} (f; z_1,  \ldots,z_{l-1}),
\end{eqnarray*}
and let $F_1 \in \cB_\infty (\tilde{Y})$
be the subspace constructed at the end
of~\ref{new_mod_2} for which  (\ref{114}) is satisfied.

Repeat the procedure inside $F_1$ by setting
\begin{eqnarray*}
  b_2 (z_1,  \ldots,z_{l-2})
&=& \bt_{z_{l-1},F_1}^{st} (b_1; z_1,  \ldots,z_{l-2}) \\
  d_2 (z_1,  \ldots,z_{l-2})
&=& \de_{z_{l-1},F_1}^{st} (d_1; z_1,  \ldots,z_{l-2}),
\end{eqnarray*}
and let $F_2 \in \cB_\infty (F_1)$
be the corresponding subspace.

Proceed by an obvious induction to get
functions $b_i$ and $d_i$ for $i = 1, \ldots, l$
and subspaces $\tilde{Y} \supset F_1 \supset \ldots \supset F_l$.
Set
\begin{eqnarray}
\bt = b_l &=&  \bt_{z_{1},F_{l-1}}^{st} (b_{l-1}) \nonumber\\
     &=&  \bt_{z_{1},F_{l-1}}^{st}
         \left(  \bt_{z_{2},F_{l-2}}^{st}
         \left (\ldots \left (  \bt_{z_{l},\tilde{Y}}^{st}
             (f; z_1,  \ldots,z_{l-1})
       \right)\ldots \right ) \right), \nonumber\\
\de = d_l &=&  \de_{z_{1},F_{l-1}}^{st} (d_{l-1}) \nonumber\\
     &=&  \de_{z_{1},F_{l-1}}^{st}
         \left(  \de_{z_{2},F_{l-2}}^{st}
         \left (\ldots \left (  \de_{z_{l},\tilde{Y}}^{st}
             (f; z_1,  \ldots,z_{l-1})
       \right)\ldots \right ) \right).
  \label{115}
\end{eqnarray}

It is easy to see, using (\ref{114}),
that with these definitions of  $\bt$ and $\de$,
the interval $[\bt, \de]$ satisfies conditions (i) and (ii)
of~\ref{def} for the  subspace $Y = F_l$.

\subsection{}
\label{stab_sub}
Using the definition of the $\bt$- and $\de$- averages
it is easy to see that
condition (i) of the definition of the tilda-spectrum
in~\ref{def} is equivalent to the following:
\begin{description}
\item[(i')]
$ \forall \theta>0 \ $
$\exists  {Y_1} \in \cB^t (Y )\ $
$\forall z_1 \in S( {Y_1})\ $
$\exists  {Y_2} \in \cB^t ( {Y_1}, z_1 ) $
$\forall z_2 \in S( {Y_2}) \ $\newline
$\exists  {Y_3} \in \cB^t ( {Y_2}, z_2)$
$\forall z_3 \in S( {Y_3}) \ $
$\ldots\ $
$\exists  {Y_l} \in \cB^t ( {Y_{l-1}}, z_{l-1})$
$\forall z_l \in S({Y_l}) \ $\newline
$  f(z_1, z_2, \ldots, z_l ) \in [\bt-\theta, \de+\theta]$.
\end{description}

We will show below that this condition  implies in fact a
stronger property, that is, the existence of a stabilizing
subspace for $f$.
Given an interval
$[\bt, \de]$  satisfying condition (i)
on a subspace $Y$
we can construct, for any  $\vt >0$,
a subspace
${G} \in \cB_\infty ({Y})$ such
that for all normalized blocks $z_1 < z_2< \ldots < z_l$ in $G$
we have
\begin{equation}
 f(z_1, z_2, \ldots, z_l ) \in
[\bt- \vt, \de+ \vt].
  \label{116}
\end{equation}

Fix $\vt' >0$ and $\eta >0$  to be defined later.
Let  $E_1 = Y_1\in \cB^t (Y)$ be
the subspace  satisfying condition (i') for $\vt'$.
Pick an arbitrary vector $u_1 \in S(E_1)$,
and let $E_2 = Y_2 \in \cB^t(Y_1, u_1)$ be the subspace
from condition (i') (again for $\vt'$).
Pick an arbitrary vector $u_2 \in S(E_2)$.

In the next step we would like to find a subspace
${E}_3 \in  \cB^t (E_2, u_2)$
which would satisfy condition (i')
in several ways: it could be taken as $Y_3$,
for vectors $z_1 = u_1$ and $z_2 = u_2$,
and it could be taken as $Y_2$,
for an   arbitrary vector $z_1$
running over some
finite $\eta$-net $\cal N$
(in the original norm)
on the sphere $S(\spn[u_1, u_2])$.
Since subspaces appearing in (i') are always of finite
codimension, it is clear that a required subspace $E_3$
exists. Then pick an arbitrary
$u_3 \in S(E_3)$.

Continuing in an obvious manner we construct
a subspace $G = \spn [u_i]$ such that
$ f(z_1, z_2, \ldots, z_l ) \in
[\bt- \vt', \de+ \vt']$
for all $z_1 < \ldots< z_l$
with $z_i$s running over all
finite $\eta$-nets
on the spheres $S(\spn[z_1,\ldots, z_k])$
(with $k = 2, 3, \ldots$).
Choosing suitable $\eta >0$ and $\vt'>0$
depending on $\vt$, we complete the proof of (\ref{116}).

\medskip
\rem
Given a sequence  $\vt_n \downarrow 0$
we can repeat the above construction for every $n$
and then pass to a diagonal subspace. We then obtain
a subspace $Z\in \cB_\infty (Y)$
with a (block) basis $\{v_i\}$
such that for every $n$ and for arbitrary
normalized blocks $n < z_1 < \ldots < z_l$
of  $\{v_i\}$ we have
$$
f(z_1, \ldots, z_l) \in [\bt - \vt_n,\de + \vt_n].
$$

\subsection{}
\label{asymp_sets}
The role of
condition (ii) of the definition of the tilda-spectrum
is to ensure the existence of large
sets of vectors on which the value
of the function $f$ is close to extremal.
In fact,  these sets  turn out
to be asymptotic in some stabilizing subspace for $f$.

Let us start by observing that condition (ii) from~\ref{def}
is equivalent to the following:
\begin{description}
\item[(ii')]
for all $ \eta >0 $ and  all subspaces
$ H_1, H_2, \ldots, H_l \in \cB_ \infty(Y ) $
we have:\\
 $\exists w_1 \in S(H_1) \ $
 $ \exists w_2 \in S(H_2), w_2 > w_1\ $
 $\ldots\ $ 
 $\exists w_l \in S(H_l), w_l > w_{l-1} \ $ 
 $   f(w_1, w_2, \ldots, w_l) \le \bt + \eta$\\
 and\\
 $\exists v_1 \in S(H_1) \ $
 $ \exists v_2 \in S(H_2), v_2 > v_1\ $
 $\ldots\ $ 
 $\exists v_l \in S(H_l), v_l > v_{l-1} \ $ 
 $   f(v_1, v_2, \ldots, v_l) \ge \de - \eta$.
\end{description}

Let $[\bt, \de]$ be in the tilda-spectrum of $f$
and let  $Z$ be the corresponding stabilizing subspace
constructed  in Remark~\ref{stab_sub}.
Condition (ii') leads to the natural definition
of sets asymptotic in $Z$.

With a fixed $\eta >0$
define $A_1 \subset S(Z)$ by
\begin{eqnarray*}
A_1 &=& \bigl\{w_1 \in S(Z)\bigm|\,\,
\forall H_2, \ldots, H_l \in \cB_\infty (Z)\\
&& \qquad
 \exists w_2 \in S(H_2), w_2> w_1\
 \exists w_3 \in S(H_3), w_3 > w_2\ \ldots\   \\
&& \qquad
 \exists w_l \in S(H_l), w_l > w_{l-1}\ \
 f(w_1, w_2, \ldots, w_l) \le \bt + \eta \bigr\}.
\end{eqnarray*}

By (ii'), the set $A_1$ has a non-empty intersection
with every subspace $H_1 \in \cB_\infty (Z)$,
hence $A_1$ is asymptotic in $Z$.

By induction, let $1 \le k < l$, and assume that
for any fixed
$w_1 < \ldots < w_{k-1}$, with
$w_i \in A_i$ for  $i = 1, \ldots, k-1$,
the set $A_k= A_k(w_1,  \ldots,  w_{k-1})$
has been defined by the formula
\begin{eqnarray}
  \label{a_k}
A_k &=& \bigl\{w_k \in S(Z)\bigm|\,\,
\forall H_{k+1}, \ldots, H_l \in \cB_\infty (Z) \nonumber\\
&& \quad
 \exists w_{k+1} \in S(H_{k+1}),  w_{k+1}> w_k\
 \exists w_{k+2} \in S(H_{k+2}), w_{k+2} > w_{k+1}\ \ldots \nonumber   \\
&& \quad
 \exists w_l \in S(H_l), w_l > w_{l-1}\ \
 f(w_1, w_2, \ldots, w_l) \le \bt + \eta \bigr\}.
\end{eqnarray}
Moreover, assume that $A_k$ is  asymptotic in $Z$.
Then for any fixed
$w_1 < \ldots < w_{k}$, with
$w_i \in A_i$ for $i = 1, \ldots, k$,
define $A_{k+1}$ by the formula analogous
to (\ref{a_k}). It clearly
follows from the form of $A_k$ that  $A_{k+1}$ is  asymptotic
in $Z$.

Similarily, we can define sets
$U_k \subset S(Z)$ for $j = k, \ldots, l$,
which are also asymptotic in $Z$,
and  such  that if  $v_1 < \ldots < v_l$ and
$v_i \in U_i$ for $i = 1, \ldots, l$
then $f(v_1, \ldots, v_l)\ge \de -\eta$.

\subsection{}
\label{uncond}
The notion of tilda-spectrum has the following
unconditionality property.
For a given function $f$ and a
finite sequence $\ep= (\ep_1, \ep_2, \ldots)$  with $\ep_1 = \pm 1$,
$\ep_2 = \pm 1$, $\ldots$,   define
the function $f_\ep$ by
$$
f_\ep (z_1, z_2, \ldots, z_l)= f (\ep_1 z_1,\ep_2 z_2, \ldots, \ep_l z_l),
$$
for normalized blocks $z_1 < \ldots, z_l$.
Then if   $[\bt, \de]$ is in the tilda-spectrum  of $f$
and $Y$ is a  corresponding spectrum subspace, then
for all the $f_\ep$s, $[\bt, \de]$ is again
a spectrum interval with the same spectrum
subspace $Y$.
Moreover, the stabilizing subspace $Z \subset Y$
of~\ref{stab_sub}
is also preserved for all  the $f_\ep$s.
Note however, that the asymptotic sets $A_i$ and  $U_i$ described
in~\ref{asymp_sets}  are not the same.

\subsection{}
\label{common-stab}
The final  important step in our discussion of tilda-spectrum
is an observation that  a construction
of  stabilizing subspaces in~\ref{stab_sub}
can be done ``almost'' simultaneously
for any countable family of
uniformly continuous real  functions,
$f_k = f_k (z_1, \ldots, z_{l_k})$.
For a subspace $\tilde{Y} \in \cB_\infty (X)$
and a sequence $\vt_k \downarrow 0$,
there exists  ${Z}\in \cB_\infty (\tilde{Y})$ such that
for every $n = 1, 2, \ldots$,
if $L_n= \max_{k \le n} l_k$, then
for arbitrary
normalized blocks $n < z_1 < \ldots <z_{L_n}$ in $Z$
we have
$$
  f_k (z_1, \ldots, z_{l_k}) \in [\bt_k- \vt_n, \de_k+ \vt_n]
\quad\mbox{\rm for}\quad 1 \le k \le n.
$$
Here $[\bt_k, \de_k]$ is an interval in the tilda-spectrum
of $f_k$ in $\tilde{Y}$.
Moreover, all the sets  $A_i^{(k)}$ and $U_i^{(k)}$
for $i = 1, 2, \ldots$,
constructed in~\ref{asymp_sets}
for the function $f_k$, are asymptotic in $Z$.

This follows
from~\ref{stab_sub} and~\ref{asymp_sets}  by the
standard diagonal procedure. The details are left for the
reader.

\section{Spaces with bounded distortions}
We now pass to the main theorem on spaces
with bounded distortions, Theorem~\ref{general_dist}.

\subsection{}
\label{bd_notat}

By passing to an infinite-dimensional subspace of $X$
and considering a suitable renorming of $X$
we may assume, without loss of generality, that $X$ has
a monotone basis.
The proof of the theorem relies on stabilization properties of
a  family of real functions
which we introduce now and fix throughout
the rest of the  argument. This family is
indexed by the set  $\QQNP$ of all   finite
sequences
with positive rational coordinates;
for $a= (a_1, \ldots, a_l) \in \QQNP$ define
the function $f_a$  on a sequence $z_1 < \ldots < z_l$
of normalized blocks by
$$
f_a (z_1, \ldots, z_l) = \Bigl\| \sum_{i=1}^l a_i z_i\Bigr\|.
$$

Let  $\{a^{\,(k)}\}$ be an enumeration   of  $\QQNP$.
We will write  $f_k$  for $ f_{a^{(k)}}$.
Fix $\vt_k \downarrow 0$  satisfying
$\vt_k < (1/2D) \|a^{\,(k)}\|_\infty$
for $k = 1, 2, \ldots$.
Let $Z \in \cB_\infty (X)$ be the stabilizing subspace
for  all the $f_k$s,
constructed in~\ref{common-stab}.
Let  $[\bt_k,\de_k]$   denote
the corresponding spectrum intervals for  $f_k$
($k = 1, 2, \ldots$).

The major role in our approach is played by two positive
real valued functions on $\QQN$ defined via  tilda-spectrum
of the $f_k$s as follows.
With fixed $\vt_k$,  $Z$, and $[\bt_k,\de_k]$,
as above, we let,
for $a = a^{(k)} \in \QQNP$,
\begin{equation}
      g(a) = \bt_k
    \qquad \mbox{\rm and }\qquad r(a) = \de_k.
  \label{g_r}
\end{equation}
These definitions can be naturally extended to all
$\QQN$, by setting,
for $\pm a = (\pm a_1, \ldots,\pm a_l)$,
$g(\pm a) = g(a)$ and $r(\pm a) = r(a)$.
We call the functions $g$ and $r$ the {\em enveloping functions\/} of $X$.

\subsection{}
  \label{g_r_ineq}

The main part of the proof of the theorem
is contained in the following proposition
concerning the behaviour of functions $g$ and $r$ for spaces
with bounded distortions.

\begin{prop}
  Assume that  a Banach space   $X$ has bounded distortions
  and let  $\dis (X) < D$. With the notation from~\ref{bd_notat},
  either there exists $1 \le p <\infty$
  such that
$$
(1/D') \left(\sum |a_i|^p \right)^{1/p}  \le g(a) \le r(a)
\le D'  \left(\sum |a_i|^p \right)^{1/p}  \quad \mbox{\rm for\ } a \in \QQN,
$$
or
$$
(1/D')\max |a_i| \le g(a) \le  r(a) \le D'\max |a_i|
\quad \mbox{\rm for\ } a \in \QQN.
$$
where $D'= 4D$.
\end{prop}

\subsection{}
\label{pf_bnd_dist}

Assuming the truth of Proposition~\ref{g_r_ineq} let us
complete the proof of the theorem.

\medskip
\noindent{\bf Proof of Therem~\ref{general_dist}}\ \
Let  $\{x_i\}$ denote the block basis for  the stabilizing subspace  $Z$,
which has been fixed in~\ref{bd_notat}.
Assume that the conclusion of Proposition~\ref{g_r_ineq}
is  satisfied for some $1 \le p \le \infty$ (with the obvious
convention for $p = \infty$).
We will then show that   $Z$
is  \aslym-$l_p$ (or \aslym-$c_0$, if $p = \infty$).

Fix $n$ and fix $\ep= \ep(n) >0$ to be defined later.
Let
$b^{\,(1)}, \ldots, b^{\,(M)}$ be an $\ep$-net in the unit
sphere $S(l_p^n)$ of $l_p^n$ in the $l_\infty$-norm, and assume
without loss of generality
that $ b^{\,(i)} \in \QQN$
for $1 \le i \le M$. Thus
for every $i$ there is $k = k_i $ such that
$ b^{\,(i)} =  a^{\,(k)}$. Let $N = \max_{1 \le i \le M} k_i$.

Let $N < y_1< \ldots <y_n$ be arbitrary
normalized blocks of $\{x_k\}$.
By \ref{common-stab} and Proposition~\ref{g_r_ineq}
we have, for every $1 \le i \le M$,
\begin{eqnarray*}
(1/D) \bigl\|b^{\,(i)}\bigr\|_p - \vt_k &\le& g(b^{\,(i)})  - \vt_k \le
\bigl\|b_1^{\,(i)} y_1 +\cdots + b_n^{\,(i)} y_n\bigr\|\\
&\le& r(b^{\,(i)})  + \vt_k \le  D \bigl\|b^{\,(i)}\bigr\|_p + \vt_k.
\end{eqnarray*}

By the choice of $\vt_k$ this implies
$$
(1/2D)
\le
\bigl\|b_1^{\,(i)} y_1 +\cdots + b_n^{\,(i)} y_n\bigr\|
 \le 2 D
$$

An easy approximation argument shows that
if $\ep$ is sufficiently
small (it is enough to take $\ep = (4Dn)^{-1}$) then
the latter estimates imply
$$
 (1/4 D) \le \|c_1 y_1 +\cdots + c_n y_n\| \le 4 D,
$$
for any $(c_1, \ldots, c_n) \in S(l_p^n)$.
Thus $ y_1< \ldots <y_n$
are $D'$-equivalent to the unit vector basis in $l_p^n$,
as required.
\qed

\section{Inequalities for enveloping functions}

The proof of Proposition~\ref{g_r_ineq} is based on specific
properties of functions $g$ and $r$ in spaces with bounded
distortions, which will be established in this section.  In what
follows we keep the notation from~\ref{bd_notat}, and in particular,
$Z \in \cB_\infty (X)$ is the stabilizing subspace for the functions
$\{f_k\}$, constructed in~\ref{common-stab}.

\subsection{}
\label{as-uncon}
A space $Y$ with  a basis $\{y_i\}$
is said to be
{\em asymptotically unconditional}
if there exists a constant $D'$ such that
for every $n$ there exists $N= N(n)$ such that
for any normalized blocks
$N < z_1 < \ldots < z_n$
of  $\{y_i\}$
and any sequence of reals
$(c_1, \ldots, c_n)$ we have
\begin{equation}
\sup_{\ep_i = \pm 1}
\|\ep_1 c_1 z_1 + \ldots + \ep_n c_n z_n\|
\le D'
\inf_{\ep_i = \pm 1}
\|\ep_1 c_1 z_1 + \ldots + \ep_n c_n z_n\|.
\label{uncon}
\end{equation}

\begin{lemma}
Let $X$ be a Banach space with
bounded distortions.
Then it contains a subspace $Y \in \cB_\infty (X)$
which is asymptotically unconditional.
\end{lemma}
\proof
Assume that $X$ has a basis $\{x_i\}$.
Given $\alpha >0$,
we will construct a block basis $\{y_i\}$
for which
(\ref{uncon}) holds with the constant
$D' = (1+\alpha)\,\dis(X)$ and
$N(n) = n$.

We will use a common and convenient notation
that if $I$ and $J$ are intervals of positive integers
then $I < J$ means $\max_{i \in I} i < \min_{j \in J} j$.
Moreover, for $x = \sum_i t_i x_i \in X$,
we set $Ix =  \sum_{i\in I} t_i x_i$.

We may assume that
$X$ does not contain $c_0$,
otherwise the proof would be finished.
For a positive integer $n$ define the norm
$\Snorm{\cdot}_n$ on $X$ by
$$
\Snorm{x}_n =\sup \Bigl\| \sum_{i=1}^n \ep_i I_i x\Bigr\|,
$$
where the supremum is taken over all intervals
$I_1 < \ldots < I_n$ and all $\ep_i = \pm 1$,
$i = 1, \ldots, n$.
Clearly,
$\|x\| \le \Snorm{x}_n \le n \|x\|$ for $x \in X$,
so $\Snorm{\cdot}_n$ is an equivalent norm on $X$.
We will show that the set
\begin{equation}
A_n = \Bigl\{x \in S(X) \bigm|
  \Snorm{x}_n \le (1+\alpha) \Bigr\}
\label{as_unc}
\end{equation}
is asymptotic in $X$.
Thus, by the Remark in~\ref{as_s2} and~(\ref{as_xxx}),
for every
$X' \in \cB_\infty(X)$ there is
$F_n \in \cB_\infty(X')$ such that
$\Snorm{x}_n \le (1+\alpha)D \|x\|$ for $x \in F_n$.
This leads to the inductive construction of subspaces
$X=F_0 \supset F_1 \supset \ldots \supset F_n\supset\ldots$
with $F_n \in \cB_\infty(F_{n-1})$
and
$\|x\| \le \Snorm{x}_n \le (1+\alpha) D \|x\|$ for $x \in F_n$
($n=1, 2,\ldots$).
It is easy to check that any block basis
$\{y_n\}$ such that $y_n \in S(F_n)$
satisfies (\ref{uncon}) with $D' = (1+\alpha)D$.

To show that the set $A_n$  given by (\ref{as_unc}) is
asymptotic  in $X$, let $W = \spn [w_i] \in \cB_\infty(X)$.
Fix $N > n$ to be defined later.
Let
$$
a_N =
\sup \Bigl\{\Bigl\|\sum_{j=1}^N \eta_j w_j\Bigr\|\  \bigm|
\eta_j = \pm 1 \Bigr\}
=\Bigl \|\sum_{j=1}^N \eta_j^0 w_j\Bigr\|.
$$
Since $X$ does not contain $c_0$, we have $a_N \to \infty$,
as $N \to \infty$.

Set
$ w = \sum_{j=1}^N \eta_j^0 w_j$.
Given intervals
$I_1 < \ldots < I_n$, let
$L_i$ be the set of all $j$
such that $\supp w_j \subset I_i$
and let
$K_i$ be the set of all $j$
such that $I_i w_j \ne 0$,
for $i=1, \ldots, n$.
For every $\ep_i = \pm 1$, with $i=1, \ldots, n$,
we have
$$
\Bigl\|\sum_{i=1}^n \ep_i I_i w \Bigr\| \le
\Bigl\|\sum_{i=1}^n \ep_i \Bigl(\sum_{j \in L_i}\eta_j^0 w_j\Bigr)\Bigr\|
 + 2n \le a_N + 2n.
$$
Therefore
$\Snorm{w / a_N}_n\le 1 + 2n / a_N$.
Thus, if $a_N >2 n/ {\alpha}$, then $w / a_N \in A_n \cap S(W)$
hence $A_n$ is asymptotic in $X$.
\qed

\subsection{}
\label{sp_func_1}
Let $X$ be a Banach space with bounded distortions,
let $Y \in \cB_\infty (X)$ be an asymptotically
unconditional subspace of $X$,
for some constant $D'$ arbitrarily
close to $ D$ which can be chosen later,
and let
$Z \in \cB_\infty (Y)$ be the
stabilizing subspace for
$\{f_k\}$ constructed in~\ref{common-stab}.

The  following lemma investigates  the behaviour
of enveloping functions
$g$ and $ r: \QQN \to \Rn{}$.

\begin{lemma}
  Assume that a Banach space $X$ has bounded distortions and let
  $\dis (X) < D$.
Then
$$
g(a) \le r(a) \le 3 D g(a)
\qquad \mbox{\rm for \ } a \in \QQN.
$$
\end{lemma}
\proof
The left hand side  inequality is obvious.
To prove the right hand side inequality,
for $i= 1, \ldots, l_k$ and $k = 1, 2, \ldots$,
let $\As{A}{i}{k}$ and $\As{U}{i}{k}$ denote the $i$th asymptotic
sets, constructed for the function $f_k$ and the subspace $Z$, as
in~\ref{asymp_sets}.

We  will prove that, with fixed $a = a^{\,(k)}  \in \QQNP$, we have
\begin{eqnarray}
  \label{sp_1}
  \exists G_1 \in \cB_\infty (Z) &&\!\!\!\!\!\!\!  \forall z_1 \in
  S(G_1)\ \exists G_2 \in \cB_\infty ( G_1,z_1)\ \forall z_2 \in
  S(G_2) \ldots \nonumber\\
      && f_k(z_1, z_2, \ldots, z_{k_l}) \le   D(\bt_k+\vt_k).
\end{eqnarray}
Applying (\ref{sp_1}) to vectors from the appropriate sets
$\As{U}{i}{k}$, for $i = 1, \ldots, l_k$, we get
$$
\de_k - \vt_k \le f_k(z_1, z_2, \ldots, z_{l_k})
\le D (\bt_k+ \vt_k).
$$

By the choice of $\vt_k$ from~\ref{bd_notat}
we have
$\vt_k \le (1/2D) \|a^{\,(k)}\|_\infty \le \bt_k/2 $,
the latter inequality yields
\begin{equation}
r(a) = \de_k \le 3 D \bt_k \le 3 D g(a).
  \label{r_av_g_av}
\end{equation}

To prove that (\ref{sp_1}) holds,
first note that since $Z$ is asymptotically
unconditional, we can assume, without
loss of generality that all  the sets
$\As{A}{i}{k}$s are symmetric
about the origin.
Consider  $\As{A}{1}{k} \subset S(Z)$.
By the Remark in~\ref{as_s2} there
exists $G_1 \in \cB_\infty (Z)$ such that
$$
(1/D)(B_X \cap G_1)\subset \mbox{\rm conv} {\As{{A}}{1}{k}}.
$$

Thus for every  $z_1 \in S(G_1)$
there exist
$w^1, \ldots, w^m \in \As{A}{1}{k}$
and $t_1, \ldots, t_m$, with $t_j >0$ and $\sum_j t_j = 1$,
such that $(1/D)z_1 = \sum_j t_j w^j$.

Now set $G_{2, 0}= G_1$ and proceed by induction
in  $j = 1, \ldots, m$. For $ j \ge 1$ consider
the set $\As{A}{2}{k}(w^j)\cap G_{2, j-1}$
constructed for the vector $w^j$.
Arguing as before
we get a  subspace
$G_{2,j} \in \cB_\infty (G_{2, j-1}, w^j)$
such that
$$
(1/D) ( B_{X} \cap G_{2,j}) \subset
\mbox{\rm conv} (\As{{A}}{2}{k}(w^j)).
$$

Let $G_2 =  G_{2,m} \in \cB_\infty (G_1, z_1)$.
For a fixed  $j = 1, \ldots, m$,
we have, for an arbitrary vector $z_2 \in S(G_2)$,
$$
(1/D)z_2 = \sum_{n=1}^{m'} s_{n,j} u^{n,j},
$$
with $u^{n,j} \in \As{A}{2}{k}(w^j)$
and $s_{n,j} >0$ for $n = 1, \ldots, m'$
and $\sum_n s_{n,j} = 1$.

We repeat the process $l_k$ times, for all subsets $\As{A}{i}{k} $
with $i = 1, 2, \ldots, l_k$.  We then have, by the definitions of
$f_k$ and of the sets $\As{A}{i}{k}$,
\begin{eqnarray*}
 f_k(z_1, z_2, \ldots, z_{l_k})
   &=& D\, \bigl\|(1/D) a_1 z_1
             +(1/D) a_2 z_2   +\cdots \bigr\|\\
      &\le& D\, \sum_j t_j \
         \bigl\| a_1 w^j +(1/D)  a_2 z_2 +\cdots \bigr\|\\
      &\le& D\,\sum_j t_j  \sum_{n} s_{n,j}\
         \bigl\| a_1 w^j +  a_2 u^{n,j} +\cdots \bigr\|\\
      &\le& \dots \le D (\bt_k +   \vt_k),
\end{eqnarray*}
which is the required estimate (\ref{sp_1}).
\qed

\subsection{}
\label{sp_func_2}
Next lemma establishes  general properties of the
function $r: \QQN \to \Rn{}$.
It says that $r$ can be
extended in a natural way to the
1-unconditional and 1-subsymmetric
norm on the space $\RRN$
of all real  finite sequences.
Moreover, denoting by $\{e_i\}$
the standard unit vector basis in  $\RRN$,
the extended norm has certain blocking
property.

Recall that a norm $|\cdot|$ on  $\RRN$
is  1-subsymmetric if
for every sequence $a = (a_1, a_2, a_3,\ldots) \in \RRN$
we have
$|(a_1,0,\ldots,0, a_2, 0,\ldots,0, a_3,\ldots)|=|a|$.

\begin{lemma}
  The function $r(\cdot)$ can be extended to
  the 1-unconditional and 1-sub\-symmet\-ric
  norm on $\RRN$.
  If   $\{u_i\}$ is a block basis of the standard
  unit vector basis
  with $r(u_i)=1$ for
  $i=1, 2, \ldots$, then for every
  $a  \in \RRN$ we have
  \begin{equation}
    g(a) \le r(d) \le r(a).
    \label{blocks}
  \end{equation}
  where for $a = (a_1, \ldots, a_l)=\sum_i a_i e_i $, by $d$ we
  denote the sequence $d = \sum_i a_i u_i$.
\end{lemma}
\proof
Consider $r(\cdot)$ as a function on $\QQN$.
It is clearly positively homogeneous and we will
prove the triangle inequality by showing that
if $a, b \in \QQN$ then
$ r(a+b) \le r(a) + r(b)$.

Indeed, let $a = a^{\,(i)} = (a_1, a_2, \ldots)$,
$b = a^{\,(j)} = (b_1, b_2,\ldots)$, and $a+b = a^{\,(m)}$.
Fix an arbitrary $\eta >0$.
Pick  arbitrary vectors $w_j \in \As{U}{j}{m} $,
for $j = 1, 2, \ldots$. Then
$$ r(a+b) - \eta \le \bigl\|(a_1 + b_1) w_1 + (a_2 + b_2) w_2 +\cdots \bigr\|,
$$
moreover,
$$
\|a_1 w_1 + a_2 w_2 +\cdots \| \le r(a)+\eta \quad \mbox{\rm and }\quad
\|b_1 w_1 + b_2 w_2 +\cdots \| \le r(b)+\eta.
$$
Since $\eta >0$ is arbitrary,
this shows the triangle inequality.

Recall that for a   sequence $a = (a_1, a_2, \ldots)$,
we set
$\pm a = (\pm a_1,\pm a_2, \ldots)$.
Then we have
$r(a) = r(\pm a)$, hence the norm $r(\cdot)$
is 1-unconditional. It is also clearly
1-subsymmetric.

To prove (\ref{blocks}), let
$u_i = \sum_{j={k_i}+1}^{k_{i+1}} b_j e_j$,
for some  $0 \le k_1 < k_2 < \ldots$,
be a  block basis with rational coefficients of the standard
unit vector basis, with $r(u_i)=1$ for
$i=1, 2, \ldots$.
Let  $a = (a_1, \ldots, a_l)=\sum_i a_i e_i $,
then  $d = \sum_i a_i u_i = (d_1, \ldots, d_{k_{l+1}})$,
where $d_j = a_i b_j$ for $ k_i < j \le k_{i+1}$ and
$i=1, 2, \ldots$.

Let $\eta >0$. There exist vectors
$w_1 < w_2 < \ldots < w_{k_{l+1}}$ in appropriate asymptotic sets
such that
$$
r(d) - \eta \le \Bigl\|\sum_{j=1}^{k_{l+1}} d_j w_j \Bigr\|.
$$
Since all the vectors belong to $Z$, we also have, from the form
of $d$,
$$
c_i = \Bigl\|  \sum_{j={k_i}+1}^{k_{i+1}} d_j w_j\Bigr\|
\le |a_i| (r(u_i) + \eta) = |a_i| (1 + \eta),
$$
for $i = 1, \ldots, l$.
By the triangle inequality and
the unconditionality of the norm $r(\cdot)$ this implies
$ r(c_1, \ldots, c_l) \le r(a) (1 + \eta)$.

Setting
$w^i = (1/c_i)  \sum_{j={k_i}+1}^{k_{i+1}} d_j w_j$
we get
$$
r(d) - \eta \le \Bigl\|\sum_{i=1}^l c_i w^i\Bigr\|
\le r(c_1, \ldots, c_l)+ \eta
$$
which combined with the previous inequality
shows the right hand side of (\ref{blocks}).

The proof of the left hand side inequality is similar.
For an arbitrary  $\eta >0$, pick vectors
$w_1' < w_2' < \ldots < w_{k_{l+1}}'$ in appropriate asymptotic sets
(for appropriate functions $f_i$)
such that for every $i = 1, \ldots, l$ one has
$$
1 - \eta = r(u_i) - \eta \le
\Bigl\| \sum_{j={k_i}+1}^{k_{i+1}} b_j w_j'\Bigr\|
\le  r(u_i) + \eta = 1 + \eta.
$$
Then
$$
\Bigl\|\sum_i a_i \sum_{j={k_i}+1}^{k_{i+1}} b_j w_j'\Bigr\|
=  \Bigl\|\sum_{j=1}^{k_{l+1}}  d_j w_j' \Bigr\|
\le r( d) + \eta.
$$

On the other hand,
setting $v_i = \sum_{j={k_i}+1}^{k_{i+1}} b_j w_j'$
for $i = 1, \ldots, l$  we have
$1 - \eta \le \|v_i\| \le 1 + \eta $.
Then
$$
\Bigl \|\sum_i  a_i \sum_{j={k_i}+1}^{k_{i+1}} b_j w_j'\Bigr\|
=  \Bigl\|\sum_{i=1}^{l} a_i v_i \Bigr\| \ge
(g (a) - \eta)(1- \eta).
$$
Combining the last two estimates  we complete the proof
of the left hand side  of (\ref{blocks}).
\qed

\subsection{}
Now we can easily complete the proof of
Proposition~\ref{g_r_ineq}.

\smallskip
\noindent{\bf Proof of  Proposition~\ref{g_r_ineq}\ \ }
Let $L$ denotes the completion of
$( \RRN, r(\cdot))$, and let $\{e_i\}$
be the standard unit vector basis in $L$.
Lemmas~\ref{sp_func_1} and~\ref{sp_func_2}
yield that all normalized block bases of $\{e_i\}$
are  (3D)-equivalent to $\{e_i\}$.
By Zippin's theorem, this implies that the
$\{e_i\}$ is equivalent to the
standard unit vector basis
in $l_p$, for some $1 \le p <\infty$
or in $c_0$. Moreover, the equivalence constant
depends on $D$ only.
\qed

\section{Asymptotic $l_p$ spaces, general properties}

We  conclude this paper  with few simple remarks on general
asymptotic $l_p$ spaces.
To avoid tiresome repetitions, when talking
about spaces $l_p$ or  asymptotic $l_p$, respectively,
we adopt the convention that the case of $p = \infty$
corresponds to the space $c_0$ or  asymptotic $c_0$,
respectively.

Let $Y$ with a basis $\{y_i\}$ be an asymptotic $l_p$ space
for $1 \le p \le \infty$, 
and let   ${\la}_p (Y)$ be the  asymptotic $l_p$ constant,
as defined in~\ref{general_dist}.

\subsection{}
\label{projections}
It is well-known and easy to see that any block subspace
of $l_p$ is complemented. 
The same is true  in  Tsirelson space and  in
its convexifications (\cf\ \eg\ [C-S]),
although in this case the argument
is much  more complicated.
An analogous fact for arbitrary \asm\ $l_p$ space
says that finite-dimensional
block subspaces far out are uniformly complemented.

More precisely, for $C > {\la}_p (Y)$,
if $ N(n)=N < z_1 < \ldots <z_n$ are normalized blocks
$C$-equivalent
to the standard unit vector basis in $l_p^n$,
then
there exists a projection $P$ from $Y$ onto $\spn[z_i]_{i=1}^n$
with $\|P\| \le 2C^2$.

Indeed, pick $z_i^* \in Y^*$ such that $\|z_i^*\|= z_i^*(z_i)=1$
for $i= 1, \ldots, n$.
Let $N < E_1 < \ldots <E_n$ be intervals of positive
integers such that $ \supp z_i\subset E_i$,
for $i= 1, \ldots, n$, and that the union of all the
$E_i$s is an interval.  For $x \in Y$  set
$$
Px = \sum_{i=1}^n  z_i^*(E_i x) z_i.
$$
Since $E_j z_i = 0$ if $i \ne j$, then $P$ is a projection.
Moreover, we have
\begin{eqnarray*}
\|Px\| &\le& C \Bigl(\sum_{i=1}^n |z_i^*(E_i x)|^p \Bigr)^{1/p}
         \le C \Bigl(\sum_{i=1}^n \|E_i x\|^p \Bigr)^{1/p}\\
 &\le& C^2 \Bigl\Vert\sum_{i=1}^n \|E_i x\|\,
       \bigl({E_i x \over \|E_i x\|}\bigr)\,\Bigr\Vert
    \le C^2 \Bigl\|\sum_{i=1}^n E_i x\Bigr\| \le 2 C^2 \|x\|,
\end{eqnarray*}
as required.

\subsection{}
\label{duality}
If $Y$ is an \asm\ $l_p$ space (for $1 < p \le \infty$) then the dual
$Y^*$ is an \asm\  $l_{p'}$ space. This follows
from~\ref{projections} by a general duality argument.
A direct calculation  is just as standard and simple
and we leave it to the reader.

\subsection{}
\label{topol}
It is easy to observe that
if $p >1$, the basis in $Y$  is shrinking
and if $p < \infty$, the basis is
boundedly complete.
Hence for $1 < p < \infty$, an
\asm\  $l_p$ space is reflexive.

Assume to the contrary that the basis $\{y_i\}$ is not shrinking.
There exists $x^* \in Y^*$ with $\|x^*\|= 1$, and $\delta >0$,
and a normalized block basis $\{u_i\}$ of  $\{y_i\}$ such that
$|x^*(u_i)| >\delta$.
 Fix $n$ to be defined later. Then for every $k$ we have
$$
\Bigl\|\sum_{i=k}^{k+n-1} u_i\Bigr\|
\ge \Bigl|x^*\Bigl(\sum_{i=k}^{k+n-1} u_i\Bigr)\Bigr|\ge n \delta.
$$
On the other hand, if $k$ is large enough, the left hand side is smaller
than or equal to $C\,n^{1/p}$. Chosing appropriate $n$ we get a
contradition, if $p >1$.

Assume the basis is not boundedly complete.
Then there exists normalized block
basis $\{u_i\}$ of  $\{y_i\}$ such that
$\sup_n \| \sum _{i=1}^n u_i\|= M <\infty$.
On the other hand if $u_k < \ldots < u_{k+n-1}$
is far enough, then
$$
\Bigl\| \sum _{i=k}^{k+n-1} u_i\Bigr\| \ge (1/C)\, n^{1/p}.
$$
If $p < \infty$,  we again come to a contradition by an appropriate choice
of $n$.

\subsection{}
\label{isomorph}
The notion of  \asm\ $l_p$ spaces is fundamentally
an isomorphic concept and it cannot be reduced to a $(1+\ep)$-
isometric one. As mentioned in~\ref{general_dist},
to the contrary to the local concept of finite
representability of $l_p$, \asm\ $l_p$ does not
imply any related  almost isometric property
of a block subspace. To be more precise
let us discuss the quantity $\la_p (Y)$ in more detail.

Note that the isometric condition
$ \la_p (Y)=1$ is
equivalent to the fact that
for every $\ep >0$ and for every $n$
there exists $N = N(\ep, n)$ such that
any normalized blocks
$N <z_1 < z_2 < \ldots< z_n$  are $(1+ \ep)$-equivalent
to the standard unit vector basis in $l_p^n$.
(In~\ref{general_dist} we called such a space
an \ai\ \as\ $l_p$ space.)

Recall that if a space $Y$ merely
satisfies a weaker condition:
for every $\ep >0$
there exists $N = N(\ep)$ such that
any two normalized blocks
$N <z_1 < z_2 $  are $(1+ \ep)$-equivalent
to the standard unit vector basis in $l_p^2$,
then,
for every $\ep >0$,
$Y$ contains  an $(1+\ep)$-isomorphic
copy of the $l_p$-space.
Let us sketch this well-known and
standard argument.

Given  $\ep >0$,
fix $\ep_i \downarrow 0$ such that
$\prod_i (1+ \ep_i) \le  (1+ \ep)$.
By an easy induction pick a sequence of normalized blocks
$u_1 < u_2 <\ldots < u_{i} <\ldots$
such that $N(\ep_i) < u_i$ for $i = 1, 2, \ldots$.
Then the block basis
$\{u_i\}$ is  $(1+ \ep)$-equivalent
to the standard unit vector basis in $l_p$.
Indeed, for any finite sequence of
scalars $\{a_i\}$ we have
\begin{eqnarray*}
  \Bigl\| \sum_{i=1}^\infty a_i u_i \Bigr\|
  &\le & (1+\ep_1)\left(|a_1|^p +
            \Bigl \| \sum_{i=2}^\infty a_i u_i \Bigr\|^p \right)^{1/p}\\
  &\le & (1+\ep_1) (1+\ep_2) \left(|a_1|^p + |a_2|^p +
            \Bigl \| \sum_{i=3}^\infty a_i u_i \Bigr\|^p \right)^{1/p}\\
   &\le& \ldots \le \prod_i (1+ \ep_i)
              \Bigl( \sum_{i=1}^\infty |a_i|^p \Bigr)^{1/p}.
\end{eqnarray*}
In a similar way we get the lower estimate, hence
$$
(1+\ep)^{-1}  \Bigl( \sum_{i=1}^\infty |a_i|^p \Bigr)^{1/p}
\le  \Bigl\| \sum_{i=1}^\infty a_i u_i \Bigr\|
\le (1+\ep)  \Bigl( \sum_{i=1}^\infty |a_i|^p \Bigr)^{1/p},
$$
as required.

Clearly, $ \la_p (Y)$ is an isomorphic invariant.
However, there exist
spaces $Y$ such that  $ \la_p (Y) < \infty$
but there is no equivalent norm
$\Snorm{\cdot}$ on $Y$ such that for
some block subspace $Z \in \cB_\infty (Y)$
the equality
$ \la_p (Y, \Snorm{\cdot} ) = 1$ would hold.
The construction above obviously yields
that this is true for every  \asm\ $l_p$
space  which does not
contain subspaces   isomorphic to $l_p$.
In particular,  Tsirelson space and
its convexifications have this property.


%
\end{document}